\documentclass[12pt,leqno,psfig,openbib]{amsart}
\usepackage{amssymb,amsmath,graphicx,amsfonts,euscript}
\usepackage{fancyhdr}
\usepackage{epsfig}
\usepackage{setspace}
\usepackage{dsfont}
\if@secthm
  \newtheorem{thm}{Theorem}[section]
\else
  \newtheorem{thm}{Theorem}[section]
\fi
\newtheorem{rem}{Remark}[section]
\newtheorem{lem}{Lemma}[section]
\newtheorem{cor}{Corollary}[section]

\numberwithin{equation}{section}

\begin{document}


\title[Very weak solutions of the Stokes problem in a convex polygon 10/12/2015]{Very weak solutions of the Stokes problem in a convex polygon}

\author[Makram Hamouda$^{1,2}$, Roger Temam$^2$ and Le Zhang$^2$ 10/12/2015]{Makram Hamouda$^{1,2}$, Roger Temam$^2$ and Le Zhang$^2$}

\address{$^1$ University of Carthage, Faculty of Sciences of Bizerte, Department of Mathematics, Zarzouna,  Tunisia.}
\email[Makram Hamouda]{mahamoud@indiana.edu}

\address{$^2$ The Institute for Scientific Computing and Applied Mathematics, Indiana University,
831 East Third Street, Bloomington, Indiana  47405, U.S.A.}
\email[Roger Temam]{temam@indiana.edu}
\email[Le Zhang]{zhangle@indiana.edu}


\baselineskip=18truept

\maketitle

\begin{abstract}Motivated by the study of the corner singularities in the so-called cavity flow, we establish in this article, the existence and uniqueness of solutions in $L^2(\Omega)^2$ for the Stokes problem in a domain $\Omega,$ when $\Omega$ is a smooth domain or a convex polygon.  We establish also a trace theorem and show that the trace of $u$ can be arbitrary in $L^2(\partial\Omega)^2.$  The results are also extended to the linear evolution Stokes problem.
\end{abstract}

\section{Introduction}\label{s1}
\par
We are interested in this article in the existence of very weak solutions for the (linear stationary) Stokes problem in a domain $\Omega$ of $\mathbb{R}^2.$  The set $\Omega$ is assumed to be bounded, regular of class $\mathcal{C}^4$, or it could be a convex polygonal domain.  More generally $\Omega$ can be what we will call a (convex) \textit{domain of polygonal type} that is $\Omega$ a piecewise $\mathcal{C}^4$ domain, for which the tangent to the boundary $\Gamma$ has a finite number of discontinuity points $S_1,\ldots , S_N,$ with a well defined left and right tangent at these points, the angle between the tangents being
$0<\alpha_j<\pi$; hence the domain $\Omega$ needs not be convex in this case, we only require the angles to be convex.

\par
Motivated by the study of a flow in such a domain $\Omega$ we start with the (stationary linear) Stokes problem, which, in its most general form reads:

\begin{equation}\label{e1.1}
\begin{cases}
-\Delta\tilde u +\nabla\tilde p=f \enspace\text{ in }\enspace\Omega,\\
\text{div }\tilde u = h\enspace\text{ in }\enspace\Omega,\\
\tilde u = g\enspace\text { on }\enspace\Gamma =\partial \Omega.
\end{cases}
\end{equation}
The emphasis in \cite{Zha15} is on the so-called lid driven cavity problem.  In this case $\Omega$ is the square $(0,1)\times (0,1), f=h=0$ and
$g=(0,0)$ at $x=0,1,$ and $y=0,$ and $g=(1,0)$, at $y=1$; the discontinuities of $g$ produce singularities and vortices at the corners $(0,1)$ and $(1,1).$  We describe this example in more detail in Section \ref{s4}.

\par
We know that if $f\in L^2(\Omega)^2, h=g=0$, then the existence and uniqueness of a solution $\tilde u\in H^1_0(\Omega)^2$ of \eqref{e1.1} is derived from the projection theorem, and $\tilde p\in L^2(\Omega)$ follows from the result of Deny-Lions \cite{DL54}; see also \cite{Tem01}, \cite{Wan93}.  If $f\in L^2(\Omega)^2$ and $h\in L^2(\Omega)$, with $\int_\Omega hdxdy =0,$ we have existence and uniqueness of $U\in H^2(\Omega)^2, P\in H^1(\Omega)$ satisfying
\begin{equation}\label{e1.2}
\begin{cases}
-\Delta U+\nabla P=f\enspace\text{ in }\enspace\Omega,\\
\text{ div }U=h\enspace\text{ in }\enspace\Omega,\\
U=0\enspace\text{ on }\enspace\Gamma.
\end{cases}
\end{equation}
In the case where $\Omega $ is smooth, this result is proved in \cite{Cat61}; see also \cite{Ghi84}.  When $\Omega$ is of polygonal type, this result is proven in \cite{Gri92}; see also \cite{Gri85}.  The difference $u=\tilde u-U, p=\tilde p-P$ is solution of the following problem which concentrates the lack of regularity on the boundary value, like for the lid driven cavity problem:
\begin{equation}\label{e1.3}
\begin{cases}
-\Delta u +\nabla p=0\text{ in }\Omega,\\
\text{ div }u=0\text{ in }\Omega,\\
u=g\text{ on }\Gamma.
\end{cases}
\end{equation}
In Section \ref{s2} we derive a trace theorem for functions $u$ in $L^2(\Omega)^2$ satisfying \eqref{e1.3}$_1$, \eqref{e1.3}$_2$, thus giving a meaning to \eqref{e1.3}$_3$.  Then in Section \ref{s3} we establish the existence and uniqueness of a solution $u\in L^2(\Omega)^2$ of \eqref{e1.3} provided
$\int_\Gamma g\cdot n\enspace d\Gamma=0$ see Theorem \ref{t3.1}.  We discuss in more details an example of lid driven cavity flow in Section \ref{s4}.  Finally in Section \ref{s5} we extend the results to the linear evolution Stokes problem.  Namely, we prove the necessary trace theorems, then, in Theorem \ref{t5.2}, we show that if $g$ is given in $L^2(0,T;L^2(\Gamma)^2)$ satisfying $\int_\Gamma g\cdot n\enspace d\Gamma=0$ for a.e. $t$ then, there exists a unique solution $u\in L^2(0,T;L^2(\Omega)^2)$ of
\begin{equation}\label{e1.3a}
\begin{cases}
&\displaystyle\frac{\partial u}{\partial t} - \Delta u +\nabla p=0\text{ in }(0,T)\times\Omega,\\
&\text{div }u=0\text{ in }\Omega,\\
&u=g\text{ on }(0,T)\times\Gamma,\\
&u(0)=0,
\end{cases}
\end{equation}
with $p\in\mathcal{D}'((0,T)\times\Omega).$

We conclude this introduction with an a priori estimate borrowed from \cite{HS15} which we reproduce for the sake of completeness.

\begin{lem}\label{l1.1}
Assume that $u,p$ and $g$ are sufficiently regular (e.g. $u\in H^2(\Omega)^2, p\in H^1(\Omega), g\in H^{3/2}(\Gamma)^2)$ and satisfy \eqref{e1.3}.  Then
\begin{equation}\label{e1.4}
|u|_{L^2(\Omega)^{2}}\leq c_1|g|_{L^2(\Gamma)^2},
\end{equation}
where the constant $c_1$ depends only on $\Omega.$
\end{lem}

\noindent\textbf{Proof}
The proof is based on a transposition argument used for general elliptic problems by Lions and Magenes \cite{LM72}.

We consider $v,q$ solution of the adjoint problem
\begin{equation}\label{e1.5}
\begin{cases}
-\Delta v+\nabla q=u \enspace\text{ in }\enspace\Omega,\\
\text{ div }v=0\enspace\text{ in }\enspace\Omega,\\
v=0\enspace\text{ on }\enspace\Gamma.
\end{cases}
\end{equation}
We know from the references quoted above that, at least, $v\in H^2(\Omega)^2, q\in H^1(\Omega)$ so that the following integrations by parts make sense.
We have, $n$ denoting the outside unit normal on $\Gamma:$
\begin{equation*}
\begin{split}
\int_\Omega|u|^2d\Omega &=\int_\Omega u\enspace (-\Delta v + \nabla q)\enspace d\Omega,\\
&=(\text{since }v=0\text{ on }\Gamma\text{ and div } u=\text{ div }v=0)\\
&=-\int_\Gamma g\displaystyle\frac{\partial v}{\partial n}\enspace d\Gamma-\int_\Omega\Delta u\enspace v\enspace d\Omega +\int_\Gamma g\cdot nq\enspace d\Gamma\\
&=\int_\Gamma(g\displaystyle\frac{\partial v}{\partial n}-g\cdot n)\enspace q\enspace d\Gamma\\
&\leq |g|_{H^{1/2}(\Gamma)^2}|\displaystyle\frac{\partial v}{\partial n}|_{H^{1/2}(\Gamma)^2} +|g\cdot n|_{H^{-1/2}(\Gamma)}|g|_{H^{1/2}(\Gamma)}\\
&\leq c|g|_{H^{1/2}(\Gamma)^2}\left(|v|_{H^2(\Omega)^2} + |q|_{H^1(\Omega)}\right)\\
&\leq (\text{using the }H^2(\Omega)\text{ regularity for }\eqref{e1.5})\\
&\leq c|g|_{H^{-1/2}(\Gamma)^2}|u|_{L^2(\Omega)}\\
&\leq c|g|_{L^2(\Gamma)^2} |u|_{L^2(\Omega)}.
\end{split}
\end{equation*}
The lemma is proven.

\begin{rem}\label{r1.1}
It follows from the proof of Lemma \ref{l1.1} that we can replace $|g|_{L^2(\Gamma)^2}$ by $|g|_{H^{-1/2}(\Gamma)^2}$ in the left-hand side of \eqref{e1.4}.  We will not use this improvement here because it is not easy to compute $|g|_{H^{-1/2}(\Gamma)^2}.$
\end{rem}

\begin{rem}\label{r1.2}
In the above we have used the $H^2-H^1$ regularity for the Stokes problem, that is $v\in H^2(\Omega)^2, q\in H^1(\Omega)$ when $u\in L^2(\Omega)^2$.
This is classical when $\Omega$ is smooth, say $\mathcal{C}^4$ (see e.g. \cite{Cat61}, \cite{Ghi84}), and is proven in \cite{Gri85}, \cite{Gri92} for the type of polygonal domains we are considering.
\end{rem}

\section{A trace theorem}\label{s2}
\par
We want to define the trace on $\Gamma$ of a function $u\in L^2(\Omega)^2$ which satisfies \eqref{e1.3}$_{1,2}$ for some distribution $p\in\mathcal{D}' (\Omega)$.  We first recall (see e.g. \cite{Tem01}), that if $u\in L^2(\Omega)^2$ and $\text{ div }u\in L^2(\Omega),$ the trace of its normal component
$\gamma_n(u)=u_n|_\Gamma$ is defined and belongs to $H^{-1/2}(\Gamma) $.  Hence we only need to define the trace of its tangential component $\gamma_\tau (u)=u\cdot\tau|_\Gamma.$

We will make use of the following result from H\'eron \cite{Her81}.

\begin{thm}\label{t2.1}
Let
\begin{equation*}
Y_2(\Omega)=\left\{v\in H^2(\Omega)^2, \text{ div }v=0\right\}.
\end{equation*}
Then a pair $(g_0,g_1)\in H^{3/2}(\Gamma)^2\times H^{1/2}(\Gamma)^2$ is the trace on $\Gamma$ of $(\gamma_0(v),\gamma_1(v)) = (v,\frac{\partial v}{\partial n})$, where $v\in Y_2(\Omega)$, if and only if
\begin{equation}\label{e1.6}
\int_\Gamma g_0\cdot n_1d\Gamma =0,
\end{equation}
and
\begin{equation}\label{e1.7}
\text{div }_\Gamma (g_0)_\tau + g_1\cdot n-2K\enspace g_0\cdot n=0,
\end{equation}
where $\text{ div }_\Gamma$ is the tangential divergence on $\Gamma, (g_0)_\tau$ the tangential component of $g_0$ and $K$ is the algebraic curvature of
$\Gamma.$
\end{thm}

When we restrict to $g_0=0,$ we obtain the following Corollary.
\begin{cor}\label{c2.1}
Let
\begin{equation*}
X_2(\Omega)=\left\{ v\in H^2(\Omega),\text{ div }v=0, v=0\text{ on }\Gamma\right\}.
\end{equation*}
Then $g_1\in H^{1/2}(\Gamma)^2$ is the trace $\gamma_1(v)$ of $\partial v/\partial n$ on $\Gamma$ where $v\in X_2(\Omega),$ if and only if
\begin{equation}\label{e1.8}
g_1\cdot n=0\quad\text{ on }\Gamma,
\end{equation}
or alternatively if and only if
\begin{equation}\label{e1.9}
g_1\in H^{1/2}_\tau(\Gamma)^2,
\end{equation}
the space of tangential vectors in $H^{1/2}(\Gamma)^2.$

Furthermore $\gamma_1$ is surjective and continuous from $X_2(\Omega)$ onto $H^{1/2}_\tau (\Gamma)^2,$ and it possesses a continuous left inverse $R$
(lifting operator), that is $R\gamma_1=I.$
\end{cor}

\textbf{Proof}
The condition \eqref{e1.8} - \eqref{e1.9} is just the restriction of \eqref{e1.6} - \eqref{e1.7} to the case $g_0=0.$ Of course the trace operator is continuous and since it is surjective, it is also continuous and one to one from the orthogonal of its kernel in $X_2(\Omega), (\text{Ker } \gamma_1)^\perp$ onto $H^{1/2}_\tau (\Gamma)^2.$ By the closed graph theorem it is bicontinuous from $(\text{Ker }\gamma_1)^\perp$ to $H^{1/2}_\tau(\Gamma)$ and its inverse $R$ is a left inverse of $\gamma_1.$ The proof is complete.
\bigskip

We now want to define the tangential trace on $\Gamma, \gamma_\tau (u)$, for a function $u$ which satisfies
\begin{equation}\label{e1.10}
\begin{cases}
\Delta u=\nabla p\enspace\text{ in }\Omega,\\
\text{div }u=0\enspace\text{ in }\Omega,
\end{cases}
\end{equation}
for some distribution $p\in\mathcal{D}'(\Omega).$  Note that, by De Rham's theory (see e.g. Deny and Lions \cite{DL54}), the set $F(\Omega)$ of $u$ in $L^2(\Omega)^2$ satisfying \eqref{e1.10} is closed in $L^2(\Omega)^2,$ and is hence a Hilbert space for the norm of $L^2(\Omega)^2.$  Note also that since $u\in L^2(\Omega)^2, \Delta u\in H^{-2}(\Omega)^2$ and by the results of Deny and Lions \cite{DL54} and the analogue of Proposition 1.2 in \cite{Tem01}, $p\in H^{-1}(\Omega).$

The construction of $\gamma_\tau (u)$ will somehow mimic the construction of $\gamma_n(u)=u\cdot n|_\Gamma$ in \cite{Tem01}.

\begin{thm}\label{t2.2}
Assume that $\Gamma$ is of class $\mathcal{C}^4$ or is of {\rm polygonal type}.  Then there exists a linear continuous operator $\gamma_\tau\in\mathcal{L} (F(\Omega), H^{-1/2}_\tau(\Gamma)^2),$ such that
\begin{equation}\label{e1.11}
\gamma_\tau u=u\cdot\tau|_\Gamma\text{ for every } u \in F(\Omega)\cap\mathcal{C}^2(\bar\Omega)^2.
\end{equation}

The following generalized Stokes formula is valid for all $u\in F(\Omega)$ and $g_1\in H^{1/2}_\tau (\Gamma)^2:$
\begin{equation}\label{e1.12}
<\gamma_\tau u,g_1> =\int_\Omega u\enspace\Delta v\enspace d\Omega,
\end{equation}
where $v$ is any function of $X_2(\Omega)$ such that $\gamma_1(v)=g_1.$

\end{thm}

\noindent\textbf{Proof}
For $u\in F(\Omega)$ and $g_1\in H^{1/2}_\tau(\Gamma)^2,$ we consider the expression
\begin{equation}\label{e1.13}
L_u(g_1)=\int_\Omega u\enspace\Delta v\enspace d\Omega,
\end{equation}
where $v$ is any function in $X_2(\Omega)$ such that $\gamma_1(v)=g_1.$  For the coherence of the definition we first need to show that the expression \eqref{e1.13} is independent of the choice of $v.$  If $v_1, v_2$ are two such choices, and $v=v_1-v_2,$ we must show that
\begin{equation}\label{e1.14}
\int_\Omega u\enspace\Delta v\enspace d\Omega =0,
\end{equation}
whenever $v\in X_2(\Omega)$ and $\gamma_1(v)=0.$  In this case
\begin{equation}\label{e2.10}
v\in\tilde X_2(\Omega) =\left\{ v\in H^2_0(\Omega),\text{ div }v=0\right\}.
\end{equation}
We prove below, in Lemma \ref{l2.1} that
\begin{equation}\label{e2.11}
\mathcal{V} =\left\{ v\in \mathcal{D}(\Omega)^2, \text{ div }v=0\right\}
\end{equation}
is dense in $\tilde X_2(\Omega).$ Then to prove \eqref{e1.14} we observe that, for $v\in\mathcal{V}$,
\begin{equation*}
\begin{split}
\int_\Omega u\enspace\Delta v\enspace d\Omega &=\enspace <u, \Delta v>_{\mathcal{D}, \mathcal{D}'} \\
&=\enspace <\Delta u, v>_{\mathcal{D}, \mathcal{D}'} = - < p,\text{ div }v>_{\mathcal{D}, \mathcal{D}',}
\end{split}
\end{equation*}
and this expression vanishes since $v\in\mathcal{V}.$

In addition, if we write \eqref{e1.13} with $v=Rg_1, R$ the left-inverse of $\gamma_1$  given by Corollary \ref{c2.1}, we observe that the expression $L_u(g_1)$ is linear continuous on $H^{1/2}_\tau (\Gamma)^2$ as indeed
\begin{equation*}
\begin{split}
|L_u(g_1)| &= |\int_\Omega u\enspace\Delta (Rg_1)\enspace d\Omega|\\
&\leq |u|_{L^2(\Omega)^2} |\Delta (Rg_1)|_{L^2(\Omega)^2}\\
&\leq c|u|_{L^2(\Omega)^2}|Rg_1|_{H^2(\Omega)^2}\\
&\leq c|g_1|_{H^{1/2}_\tau(\Gamma)^2}.
\end{split}
\end{equation*}

Finally, to prove \eqref{e1.11} we observe that if $u\in F(\Omega)\cap\mathcal{C}^2(\bar\Omega)^2$ and $v\in X_2(\Omega)$ then the following integrations by parts are legitimate:
\begin{equation*}
\begin{split}
\int_\Omega u\enspace\Delta v\enspace d\Omega &=\int_\Omega\Delta u\enspace v +\int_\Gamma(u\displaystyle\frac{\partial v}{\partial n}-\displaystyle\frac{\partial u}{\partial n}v)d\Gamma\\
&=\int_\Gamma u\displaystyle\frac{\partial v}{\partial n}d\Gamma +\int_\Omega \nabla p\enspace v\enspace d\Omega\\
&=\int_\Gamma u\enspace\displaystyle\frac{\partial v}{\partial n}\enspace d\Gamma +\int_T p\enspace v\enspace n\enspace d\Gamma-\int_\Omega p
\text{ div }v\enspace d\Omega\\
&=\int_\Gamma u\enspace\displaystyle\frac{\partial v}{\partial n}\enspace d\Gamma.
\end{split}
\end{equation*}

We are left with proving the following:

\begin{lem}\label{l2.1}
$\mathcal{V}$ is dense in $\tilde X_2(\Omega).$

\end{lem}

\noindent\textbf{Proof}
It is clear that $\mathcal{V}\subset \tilde X_2(\Omega)$.  We will prove the lemma by showing that any linear continuous form $\ell$ on $\tilde X_2(\Omega)$ which vanishes on $\mathcal{V}$ is equal to zero.

We observe that $\tilde X_2(\Omega)$ is a closed subspace of $(H^2(\Omega)\cap H^1_0(\Omega))^2$ and, noticing that $|\Delta v|_{L^2(\Omega)}$ is a norm on
$H^2(\Omega)\cap H^1_0(\Omega),$ we see that $\ell$ is necessarily of the form
\begin{equation*}
\ell (v)=\int_\Omega\phi\enspace\Delta v\enspace d\Omega,
\end{equation*}
for some $\phi\in L^2(\Omega)^2.$  We now write that $\ell$ vanishes on $\mathcal{V}$:

\begin{equation*}
\begin{split}
\ell(v) &= \int_\Omega\phi\enspace\Delta v\enspace d\Omega =\enspace <\phi,\Delta v>_{\mathcal{D},\mathcal{D}'},\\
&=\enspace <\Delta\phi, v>_{\mathcal{D},\mathcal{D}'},\quad \forall v\in\mathcal{V}.
\end{split}
\end{equation*}
But this classically implies that $\Delta\phi =\nabla\pi$ for some $\pi\in\mathcal{D}'(\Omega)$ which actually belongs to $H^{-2}(\Omega)$ as observed above.

Now for any $v\in\tilde X_2(\Omega),\nabla\pi=\Delta\phi\in H^{-2}$ and
\begin{equation*}
\ell (v)= \enspace <\nabla\pi, v>_{H^{-2}, H^2}.
\end{equation*}
This expression vanishes because $v\in H^2_0(\Omega)^2.$ Indeed using the fact that $\nabla\pi\in H^{-2}(\Omega),$ we see that
\begin{equation*}
<\nabla\pi, v>_{H^{-2},\enspace H^2_0}\enspace = - <\pi,\nabla v>_{H^{-1},\, H^1_0}\enspace =\enspace \lim_{j\rightarrow\infty} <\pi, \nabla v_j>_{H^{-1},\, H^1_0},
\end{equation*}
where $v_j\in\mathcal{D}(\Omega)^2$ converges to $v$ in $H^2_0(\Omega)^2$ (\textit{not necessarily in} $\tilde X_2(\Omega)$).  In the end
\begin{equation*}
<\nabla\pi, v>_{H^{-2},\enspace H^2_0}\enspace =\enspace <\pi,\nabla v> \enspace =0.
\end{equation*}

\section{The main theorem (Time independent case)}\label{s3}
\par

Our aim is now to prove the existence and uniqueness of a solution $u\in L^2(\Omega)^2$ for \eqref{e1.3} when $g$ is given in $L^2(\Gamma)^2.$

\begin{thm}\label{t3.1}
We assume that $\Omega$ is of class $\mathcal{C}^4$ or is of {\rm polygonal type}, and that $g$ is given in $L^2(\Gamma)^2$ satisfying
\begin{equation}\label{e1.15}
\int_\Gamma g\cdot n\enspace d\Gamma =0.
\end{equation}

Then there exists a unique solution $u\in L^2(\Omega)^2$ satisfying \eqref{e1.3} for some $p\in\mathcal{D}'(\Omega).$
\end{thm}

\noindent\textbf{Proof}
It is easy to construct a sequence $\tilde g_j\in H^{3/2}(\Gamma)^2$ (or possibly more regular), which converges to $g$ in $L^2(\Gamma)^2.$  Considering then
\begin{equation}\label{e3.1a}
g_j=\tilde g_j-\frac{1}{|\Gamma|}\left(\int_\Gamma \tilde g_j\cdot n\enspace dT\right) n,
\end{equation}
we see that the $g_j$ belong (at least) to $H^{3/2}(\Gamma)^2,$ satisfy \eqref{e1.15} and converge to $g$ in $L^2(\Gamma)^2.$

For each $j,$ thanks to \cite{Cat61} when $\Omega$ is of class $\mathcal{C}^4,$ and to \cite{Gri92} when $\Omega$ is of polygonal type, we infer the existence of $(u_j,p_j)\in H^2(\Omega)^2\times H^1(\Omega),$ satisfying \eqref{e1.3} with $g$ replaced by $g_j.$  We see, thanks to Lemma \ref{l1.1}, that the sequence $u_j$ is bounded in $L^2(\Omega)^2.$  More precisely, as observed before, the sequence $u_j$ is bounded in $F(\Omega).$  As $j\rightarrow\infty$, we infer the existence of $u\in F(\Omega)$ and a subsequence still denoted $j$ such that
\begin{equation*}
u_j\rightarrow u\text{ weakly in }F(\Omega),
\end{equation*}
that is weakly in $L^2(\Omega)^2,$ and $u$ satisfies
\begin{equation*}
\begin{split}
-\Delta u +\nabla p&=0\text{ in }\Omega,\\
\text{div }u&=0\text{ in }\Omega,
\end{split}
\end{equation*}
for some distribution $p\in\mathcal{D}'(\Omega).$

In addition, the trace theorem from \cite{Tem01} for $\gamma_n,$ and the trace Theorem \ref{t2.2} above for $\gamma_\tau,$ tell us that
$\gamma_0(u_j)=g_j$ converges to $\gamma_0(u)=g$ in $H^{-1/2}(\Gamma)^2.$  Hence $u$ satisfies \eqref{e1.3}.

It remains to show the uniqueness of solution of \eqref{e1.3}.  If $(u_1,p_1,), (u_2,p_2)$ are two solutions of \eqref{e1.3} and if $u=u_1-u_2, p=p_1-p_2,$ then
\begin{equation}\label{e1.16}
\begin{cases}
&-\Delta u+\nabla p=0\enspace\text{ in }\Omega,\\
&\text{div }u=0\enspace\text{ in }\Omega,\\
&\gamma_0(u)=0\enspace\text{ on }\Gamma.
\end{cases}
\end{equation}
We must show that $u=p=0.$  We consider $(v,q)$ defined by \eqref{e1.5} as in the proof of Lemma \ref{l1.1}.  Note that $v\in X_2(\Omega),$ and $q\in H^1(\Omega).$  We then write
\begin{equation}\label{e3.2a}
\begin{split}
\int_\Omega |u|^2d\Omega &=-\int_\Omega u\enspace\Delta u\enspace d\Omega +\int_\Omega u\enspace\nabla q\enspace d\Omega\\
&= (\text{using \eqref{e1.12} and }\gamma_0(u)=0)\\
&= \int_\Omega u\enspace\nabla q\enspace d\Omega.
\end{split}
\end{equation}

According to the integration by parts formula I. (1.9) in \cite{Tem01}, this last expression is equal to
\begin{equation*}
-(\text{ div }u,\varphi) + (\gamma_n(u),\gamma_0(q)),
\end{equation*}
and, it thus vanishes since $\gamma_n(u) =0$ by \eqref{e1.16}$_3$ and $\text{ div }u=0$ by \eqref{e1.16}$_2$.

We conclude that $u=0$, thus proving the uniqueness.

\begin{rem}\label{r3.1}
Another way to approach the problem \eqref{e1.3} would be to introduce the stream function $\Psi$ such that $u=\partial\Psi/\partial y, v= - \partial\Psi/\partial x$ and then \eqref{e1.3} reduces to a biharmonic problem for $\Psi$

\begin{equation}\label{e1.17}
\begin{cases}
&\Delta^2\Psi =0\text{ in }\Omega,\\
&\Psi=0\enspace , \enspace\displaystyle\frac{\partial\Psi}{\partial n} = g\text{ on }\Gamma.
\end{cases}
\end{equation}
We would then be looking for a solution $\Psi\in H^1(\Omega)$ of problem\eqref{e1.17} properly formulated.  When $\Omega$ is smooth such problems are treated in Lions and Magenes \cite{LM72}, although the problem \eqref{e1.17} is not explicitly mentioned in \cite{LM72}.  When $\Omega$ is a convex polygon or a domain of \textit{polygonal type}, the methods of \cite{Gri85}, \cite{Gri92} might apply but this remains to be done.
\end{rem}

\begin{rem}\label{r3.2}
Various results describing the behavior of a fluid in a domain with corners appear in \cite{BDT13}, \cite{DT15}, \cite{Ser83}.  In the engineering and fluid mechanics literature see \cite{Kel83}, \cite{Mof63}, \cite{SLB89} and \cite{SSB91}.
\end{rem}

\section{Example for a lid driven cavity flow}\label{s4}
We are interested in the cavity flow where $\Omega = (0,1) \times (0,1)$ and the velocity on $\Gamma$ is $(0,0)$ at $x=0,1,$ and $y=0$ and, at $y=1, g=(1,0).$
This is a classical model problem in computational fluid dynamics which has been the object of many studies see e.g. \cite{BD73},  \cite{Bot12}, \cite{Bur66}, \cite{ECG05}, \cite{ESW05}, \cite{FPT71a}, \cite{FPT71b}, \cite{Gar11}, \cite{GGS82}, \cite{GMN81}, \cite{GMPQ02} - in dimension 3 -, \cite{She89} and the references therein.  See also \cite{Cuv86}, \cite{Dav06}, \cite{Dun85}.  The singularities at the corners $(0,1), (1,1)$ remain a substantial computational difficulty.  In \cite{Zha15} the author addresses this difficulty by replacing $g$ by a continuous function $g_\varepsilon$ which converges to $g$ in $L^2(\Gamma)^2$ as $\varepsilon\rightarrow 0.$  Such an approach has been successfully applied to the Korteweg de Vries and nonlinear Schr\"oedinger equations, to deal with incompatible data; see \cite{FQT12} and \cite{QT12}; see also \cite{CQT11a}, \cite{CQT11b} and see \cite{Tem01} regarding incompatible initial data.

We can approximate $g$ by $g_\varepsilon$ which is identical to $g$ except for the first component which is equal to
\begin{equation}\label{e4.2}
1-\sigma (x)e^{-x/\varepsilon} -\sigma (1-x)e^{-(1-x)/\varepsilon},
\end{equation}
where $\sigma$ is a smooth function
\begin{equation*}
\sigma(x)=
\begin{cases}
&1\enspace ,\enspace 0\leq x\leq 1/2\\
&\in [0,1],\enspace \frac{1}{2}\leq x\leq \frac{3}{4}\\
&0\enspace ,\enspace \frac{3}{4}\leq x\leq 1.
\end{cases}
\end{equation*}
Both $g$ and $g_\varepsilon$ satisfy the necessary conditions
\begin{equation}\label{e4.3}
\int_\Gamma g\cdot n\enspace d\Gamma =\int_\Gamma g_\varepsilon\cdot n\enspace d\Gamma =0.
\end{equation}
It is clear that $g_\varepsilon$ converges to $g$ in $L^2(\Gamma)^2$ as $\varepsilon\rightarrow 0.$  In view of Theorem \ref{t3.1}, there exists a unique solution $(u,p)$ to the problem \eqref{e1.3} with $u\in L^2(\Omega)^2$, and a (unique) solution $(u_\varepsilon, p_\varepsilon)$ to the problem \eqref{e1.3} with $g$ replaced by $g_\varepsilon$, and $u_\varepsilon\in L^2(\Omega)^2$; as usual the uniqueness of $p_1,p_\varepsilon$ is meant up to the addition of a constant. In addition, when $\varepsilon\rightarrow 0,$
\begin{equation}\label{e4.4}
g_\varepsilon\rightarrow g\text{ in }L^2(\Gamma)^2,
\end{equation}
and consequently
\begin{equation}\label{e4.5}
u_\varepsilon\rightarrow u\text{ in }L^2(\Omega)^2.
\end{equation}

\begin{rem}\label{r4.1}
If one wants to focus on one of the corner singularities only, one can consider the following variations of the cavity problem

\begin{itemize}
\item[i)] For the corner $(0,1): g^\flat$ is equal to $g$ except the first component which is equal to $y$ at $x=1$.  Then $g^\flat_\varepsilon$ is the same as $g^\flat$ except the first component which is equal to
    \begin{equation*}
    1-\sigma (y)e^{y/\varepsilon},
    \end{equation*}
    on $y=1.$
\item[ii)] For the corners $(1,1):g^\sharp$ is equal to $g$ except the first component which is equal to $y$ at $x=0.$  Then $g^\sharp_\varepsilon$ is the same as $g^\sharp$ except the first component which is equal to
    \begin{equation*}
    1-\sigma(1-y)e^{-(1-y)/\varepsilon},
    \end{equation*}
    on $y=1.$
\end{itemize}

It is clear that the analogues of \eqref{e4.3}, \eqref{e4.4}, \eqref{e4.5} are still valid in this case.
\end{rem}

\section{The Time dependant case}\label{s5}
\par

We now want to derive the analogue of Theorem \ref{r3.1} in the time dependant case.  We consider $T>0$ and set $Q_T=\Omega\times (0,T), \Gamma_T=\Gamma x(0,T),$ and we are interested in very weak solutions (in $L^2(Q_T)$) of the linearized evaluation Stokes problem
\begin{equation}\label{e5.1}
\begin{cases}
&\displaystyle\frac{\partial\tilde{u}}{\partial t} - \Delta\tilde{u} +\nabla\tilde{p}=f \enspace\text{ in } Q_T,\\
&\text{div } \tilde{u}=h\enspace\text{ in }Q_T,\\
&\tilde u=g\enspace\text{ on }\Gamma_T,\\
&\tilde u(0)=u_0\enspace\text{ in }\Omega.
\end{cases}
\end{equation}
As in the stationary case, thanks to the result of \cite{Sol64a}, \cite{Sol64b}, we can introduce a lifting of $f,h,$ and $u_0,$ by considering the solution
$U,P$ of
\begin{equation}\label{e5.2}
\begin{cases}
&\displaystyle\frac{\partial U}{\partial t} -\Delta U+\nabla P=f\enspace\text{ in }Q_T,\\
&\text{div } U=h\enspace\text{ in }Q_T,\\
&U=0\enspace\text{ on }\Gamma_T,\\
&U(0)=u_0.
\end{cases}
\end{equation}

The results of \cite{Sol64a}, \cite{Sol64b}, guarantee enough regularity for $f,h, u_0$\footnote{As we said, motivated by the lid driven cavity flow, we are interested in low regularity in $g$ but assume $f,h,u_0$ as smooth as desirable}, and then $u=\tilde{u}-U, p=\tilde{p}-P$ are solutions of
\begin{equation}\label{e5.3}
\begin{cases}
&\displaystyle\frac{\partial u}{\partial t}-\Delta u+\nabla p=0\enspace\text{ in }Q_T,\\
&\text{div } u=0\enspace\text{ in }Q_T,\\
&u=g\enspace\text{ on }\Gamma_T,\\
&u(0)=0.
\end{cases}
\end{equation}
For suitable $g's,$ we are interested in very weak solutions of \eqref{e5.3} of the form $u\in L^2(Q_T)$.

Firstly we must show that \eqref{e5.3}$_4$ makes sense when $u\in L^2(Q_T)$ satisfies \eqref{e5.3}$_{1-3}$ and, say, $g\in L^2(\Gamma_T).$  Let $V_2$ be the closure of $\mathcal{V}$ in
$H^2_0(\Omega)^2,$ then $V_2\subset H\subset V^\prime_2$ with continuous injections and each space dense in the next one and, $\forall v\in L^2(0,T;\mathcal{V}),$
\begin{equation}\label{e5.4}
\displaystyle\frac{d}{dt}(u,v)=<\Delta u,v>-<\nabla p,v>=<u,\Delta v>,
\end{equation}
and by continuity, \eqref{e5.4} holds for every $v\in L^2(0,T;V_2)$; and we conclude that
\begin{equation}\label{e5.5}
\displaystyle\frac{\partial u}{\partial t}\in L^2(0,T;V^\prime_2)
\end{equation}
and $u(0)$ makes sense, since $u$ is a.e. equal to a continuous function from $[0,T]$ into $V^\prime_2$.
\vskip0.2in

\noindent\textbf{A Priori estimate}\goodbreak

Our first result now will be an analogue of Lemma \ref{l2.1}.

\begin{lem}\label{l5.1}
Assume that $u,p,g$ and $u_0$ are sufficiently regular (e.g. $u\in L^2(0,T;H^2 (\Omega)^2, p\in L^2(0,T;H^1(\Omega)), u_0\in H^1(\Omega)^2)$), and satisfy
\eqref{e5.3}.  Then
\begin{equation}\label{e5.5}
|u|_{L^2(0,T;L^2(\Omega)^2)} \leq c_2(|u_0|_{L^2(\Omega)^2} + |g|_{L^2(0,T;L^2(\Omega)^2)})
\end{equation}
where the constant $c_2$ depends only on $\Omega$ and $T.$
\end{lem}

\noindent\textbf{Proof}

The proof is based again on a transposition argument as in Lemma \ref{l1.1}.

We consider $v,q$ solutions of the adjoint system
\begin{equation}\label{e5.7}
\begin{cases}
&-\displaystyle\frac{\partial v}{\partial t} - \Delta v+\nabla q=u\enspace\text{ in }Q_T,\\
&\text{div }u=0\enspace\text{ in }Q_T,\\
&v=0\enspace\text{ on }\Gamma_T,\\
&v(T)=0,
\end{cases}
\end{equation}
and we write
\begin{equation*}
\begin{split}
\int_{Q_{T}}u^2dxdt=&\int_{Q_{T}} u\enspace (-\displaystyle\frac{\partial v}{\partial t} -\Delta v+\nabla q)\enspace dxdt\\
=&\int_\Omega u_0\enspace v(0)dx + \int_{Q_{T}}(\displaystyle\frac{\partial u}{\partial t}-\Delta u)\enspace vdxdt
-\int_{\Gamma_{T}}g\enspace (\displaystyle\frac{\partial v}{\partial n}-nq)\enspace d\Gamma_T\\
=&\int_\Gamma u_0\enspace v(0)dx -\int_{Q_{T}}\nabla p\enspace v\enspace d Q_T - \int_{\Gamma_{T}}g\enspace (\displaystyle\frac{\partial v}{\partial n}-nq)\enspace d\Gamma_T\\
=&\int_\Omega u_0\enspace v(0)dx -\int_{\Gamma_{T}}g\enspace (\displaystyle\frac{\partial v}{\partial n}-nq)\enspace d\Gamma_T\\
&\leq |u_0|_{L^2(\Omega)^2}|v(0)|_{L^2(\Omega)^2}\\
&+ |g|_{L^2(0,T;L^2(\Gamma_T)^2)} |\displaystyle\frac{\partial v}{\partial n}-nq|_{L^2(0,T;L^2(\Gamma_T)^2)}.
\end{split}
\end{equation*}

According to the regularity results of \cite{Sol64a}, \cite{Sol64b},
\begin{equation}\label{e5.8}
|v|_{L^2(0,T;H^2(\Omega)^2)} + |q|_{L^2(0,T;H^1(\Omega))} +|\frac{\partial v}{\partial t}|\enspace_{L^2(0,T;L^2(\Omega)^2)}
\leq c |u|_{L^2(0,T;L^2(\Omega)^2)},
\end{equation}
for some constant $c$ depending on $\Omega$ and $T,$ and \eqref{e5.5} follows.

\begin{rem}\label{r4.1}
In relation with Remark \ref{r1.2}\ :   the result of regularity \eqref{e5.8} is proved in \cite{Sol64a}, \cite{Sol64b} when $\Omega$ is smooth.  Alternatively, and to cover the case $\Omega$ of polygonal type, we can write equation \eqref{e5.7} (linear Stokes evolution problem with homogeneous boundary conditions), in the abstract form
\begin{equation}\label{e5.9}
-\displaystyle\frac{dv}{dt} + Av=u,\enspace 0<t<T,\enspace v(T)=0,
\end{equation}
where $A$ is the abstract Stokes problem in the space $H=\left\{ v\in L^2(\Omega)^2,\text{ div }v=0, v\cdot n=0\right.$ on $\left.\Gamma\right\}.$  By abstract (functional analysis) argument, we see that
\begin{equation*}
\displaystyle\frac{dv}{dt},\enspace Av\in L^2(0,T;H),
\end{equation*}
and $Av\in L^2(0,T;H)$ gives $v\in L^2(0,T;H^2(\Omega)^2), q\in L^2(0,T;H^1(\Omega)),$ by invoking the stationary regularity results recalled in Remark \ref{r1.2}.  A result similar to \eqref{e5.5} is valid for the full system \eqref{e5.1}.
\end{rem}

\noindent\textbf{The trace issues}

We next deal with the trace issues.

In relation with \eqref{e5.3} we want to consider the set $\mathcal{F}(Q_T)$ of functions $u$ in $L^2(Q_T)^2$ such that, for some distribution $p$ on $Q_T, p\in \mathcal{D}'(Q_T)$, we have
\begin{equation}\label{e5.10}
\begin{cases}
&\displaystyle\frac{\partial u}{\partial t}-\Delta u+\nabla p=0\enspace\text{ in } Q_T,\\
&\text{div }u=0\enspace\text{ in } Q_T\\
&u(0)=0\enspace\text{ in }\Omega.
\end{cases}
\end{equation}
Recall that we have shown in \eqref{e5.3}-\eqref{e5.5} that if $u\in L^2(Q_T)^2$ satisfies \eqref{e5.10}$_{1-2}$ then $u$ is a.e. equal to a continuous function from $[0,T]$ into $V'_2.$  Hence \eqref{e5.10}$_3$ makes sense and the definition of $\mathcal{F}(Q_T)$ is valid.  Furthermore using the current theory \cite{Rha84} as in the stationary case, we see that $\mathcal{F}(Q_T)$ is closed in $Q_T.$  Indeed if \eqref{e5.10} holds for a sequence
$u_n\in L^2(Q_T)^2, p_n\in\mathcal{D}'(Q_T)$, and $u_n\rightarrow u$ in $L^2(Q_T)^2,$ then $\nabla p_n$ converges to $F=\Delta u-\partial u/\partial t$ in $\mathcal{D}'(Q_T)$ and necessarily $F=\nabla p$ for some $p\in \mathcal{D}'(Q_T),$ so that \eqref{e5.10} holds for $u$ and $p.$  Similarly $\text{div }u_n=0$ converges to $\text{div }u$ in $\mathcal{D}'(Q_T)$ so that $\text{div }u=0$ and the trace $u_n(0)$ being continuous, it converges to $u(0)$ in $V'_2.$  Hence \eqref{e5.10} holds and $u\in\mathcal{F}(Q_T)$.

Concerning the trace on $\Gamma_T$ of a function $u$ in $\mathcal{F}(Q_T)$, we first observe that the normal component of $u,\gamma_n(u)=u\cdot n|_{\Gamma_{T}}$ is defined and belongs to $L^2(0,T;H^{-1/2}(\Gamma_T)^2),$ according to a standard trace result in the Navier-Stokes theory \cite{Tem01}.  So the issue is to define the tangential component of $u$ on $\Gamma_T,\gamma_\tau (u).$  As in the stationary case, we will define $\gamma_\tau(u)$ by its duality with the trace on $\Gamma_T$ of a suitable function $v.$  More precisely consider
\begin{equation*}
\mathcal{X}_2(Q_T) = \left\{ v\in L^2(0,T;H^2(\Omega)^2,\text{ div }v=0, v=0\text{ on }\Gamma_{T}\right\},
\end{equation*}
and the subspace
\begin{equation*}
\begin{split}
\mathcal{Y}_2(Q_T) = \{ &v\in L^2(0,T;H^2(\Omega)^2),\text{ div }v=0,\\
&\displaystyle\frac{\partial v}{\partial t}\in L^2(0,T;L^2(\Omega)^2), v=0\text{ on }\Gamma_T,v(T)=0\}.
\end{split}
\end{equation*}
Using the lifting operator $R$ from Corollary \ref{c2.1}, we see that $\gamma_1$ is a surjective operator from $\mathcal{X}_2(Q_T)$ onto $L^2(0,T;H^{1/2}_\tau (\Gamma)^2).$
Also it is elementary to see that $\mathcal{Y}_2(Q_T)$ is dense in $\mathcal{X}_2(Q_T)$ equipped with the norm $|v|_{L^2(0,T;H^2(\Omega)^2)}$ so that the traces $\gamma_1(v)$ for $v\in\mathcal{Y}_2(Q_T)$ are dense in $L^2(0,T; H^{1/2}_\tau(\Gamma)^2).$

Now for $g_1\in\gamma_1(\mathcal{Y}_2(Q_T))$ let $v$ be one of the functions in $\mathcal{Y}_2(Q_T)$ such that $\gamma_1v=\partial v/\partial n|_{\Gamma_{T}}=g_1.$  For $u$ satisfying \eqref{e5.10}, consider the expression
\begin{equation}\label{e5.11}
\mathcal{L}_u(g_{1})=-\int_{Q_{T}} u\enspace(\frac{\partial v}{\partial t}+\Delta v)\enspace dxdt.
\end{equation}
For $u$ and $v$ smooth we have by integration by parts, and using Green's formula,
\begin{equation*}
\begin{split}
\mathcal{L}_u(g_1) &= \int_{Q_{T}}v\enspace (\frac{\partial u}{\partial t}-\Delta u)\enspace dxdt +\int_{\Gamma_{T}}u
\enspace\frac{\partial v}{\partial n}\enspace d\Gamma_T\\
&= - \int_{Q_{T}}v\enspace\nabla p\enspace dx dt +\int_{\Gamma_{T}} u\enspace\frac{\partial v}{\partial n}\enspace d\Gamma_T\\
&=\int_{\Gamma_{T}}u\enspace\frac{\partial v}{\partial n}\enspace d\Gamma_T=\int_{\Gamma_{T}}u\enspace g_1\enspace d\Gamma_T.
\end{split}
\end{equation*}
Hence the expression \eqref{e5.11} has the potential to define and characterize the tangential components of $u$ in $L^2(0,T;H^{-1/2}_\tau(\Gamma_{\tau}))$ since the $g$ under consideration are dense in $L^2(0,T;H^{1/2}_\tau(\Gamma_T)^2).$

Our next task is to show that the expression $\mathcal{L}_u(g_1)$ is independent of the choice of $v\in{\mathcal{Y}}_2(Q_T)$ such that $\frac{\partial v}{\partial n} = g$ on $\Gamma_T.$  Consider two such functions $v_1,v_2$ and their difference $v=v_1-v_2$.  We must show that
\begin{equation}\label{e5.12}
\int_{Q_{T}}u\enspace (\displaystyle\frac{\partial v}{\partial t}+\Delta v)\enspace dxdt=0,
\end{equation}
when $u\in\mathcal{F}(Q_{T})$ and $v\in\mathcal{Y}_2(Q_{T})$ satisfies
$\partial v/\partial n=0$ on $\Gamma_T.$  Then such a $v$ belongs to $L^2(0,T;H^2_0(\Omega)^2)$, that is $L^2(0,T;\tilde{X}_2(\Omega)),\tilde{X}_2(\Omega)$ as in \eqref{e2.10}.  For $v$ in $L^2(0,T;\mathcal{V}),$ the expression \eqref{e5.12} vanishes because it is equal to
\begin{equation*}
\begin{split}
<&-\displaystyle\frac{\partial u}{\partial t}+\Delta u,v>_{\mathcal{D}'(Q_T),\mathcal{D}(Q_T)}\\
&=\enspace <-\nabla p,v>_{\mathcal{D}'(Q_T),\mathcal{D}(Q_T)}\enspace =\enspace 0.
\end{split}
\end{equation*}
Since we showed in Lemma \ref{l2.1} that $\mathcal{V}$ is dense in $\tilde{X}_2(\Omega),$ we infer that $L^2(0,T;\mathcal{V})$ is dense in $L^2(0,T;\tilde{X}_2(\Omega))$ and \eqref{e5.12} holds for any $v$ in $\mathcal{Y}_2(Q_T)$ satisfying $\partial v/\partial n=0$ on $\Gamma_T.$

We then have the analogue of Theorem \ref{t2.2}.

\begin{thm}\label{t5.1}
Under the hypotheses of Theorem \ref{t2.1} $(\Gamma$ of class $\mathcal{C}^4$ or of {\rm polygonal type}), there exists a linear continuous operator
$\gamma_\tau\in\mathcal{L}(\mathcal{F}(Q_T), L^2(0,T;H^{-1/2}_\tau(\Gamma)^2))$ such that
\begin{equation*}
\gamma_\tau u=u\cdot\tau|_{\Gamma_{T}}\text{ for every }u\in\mathcal{F}(Q_T)\cap\mathcal{C}^2(\bar{Q}_T).
\end{equation*}
\end{thm}
\textit{The following generalized Stokes formula is valid for all} $u\in\mathcal{F}(Q_T)$ \textit{and} $g_1\in\gamma_1(\mathcal{Y}_2(Q_T))\subset L^2(0,T;H^{-1/2}_\tau(\Gamma)^2))$:
\begin{equation}\label{e5.13}
<\gamma_\tau u,g_1> = \int_{Q_{T}}u\enspace (\frac{\partial v}{\partial t}+\Delta v)\enspace dxdt,
\end{equation}
where $v$ is any function of $\mathcal{Y}_2(Q_T)$ such that $\gamma_1(v)=g_1.$

All statements have been proven or are proven as in Section \ref{s2}.

\vskip0.3in

\noindent\textbf{The existence and uniqueness theorem}

\begin{thm}\label{t5.2}
We assume that $\Omega$ is of class $\mathcal{C}^4$ or is of {\rm polygonal type} and that $g$ is given in $L^2(0,T;L^2(\Gamma)^2)$ satisfying \eqref{e1.15} for a.e. $t\in[0,T].$

Then there exists a unique function $u\in L^2(Q_T)^2$ satisfying \eqref{e5.3} for some $p\in\mathcal{D}'(Q_T).$
\end{thm}

\noindent\textbf{Proof}

We approach $g$ by a sequence $g_j\in L^2(0,T;H^{3/2}(\Gamma)^2)$ as in \eqref{e3.1a}, where $g_j$ converges to $g$ in $L^2(0,T;H^{-1/2}(\Gamma)^2)$ as $j\rightarrow\infty.$ For each $j$ we find, by \cite{Sol64a}, \cite{Sol64b}, $u_j,p_j\in L^2(0,T;H^2(\Omega)^2)\times L^2(0,T;H^1(\Omega))$ satisfying
\eqref{e5.3}.  The estimates provided by Lemma \ref{l5.1} show that the sequence $u_j$ is bounded in $L^2(Q_T)^2.$  Therefore $u_j$ contains a subsequence weakly convergent in $L^2(Q_T)^2$ to some limit $u,$ and $u$ satisfies \eqref{e5.3} for some $p\in\mathcal{D}'(Q_T),$ thanks to \cite{Rha84}.

We are left with the uniqueness, that is proving that $u=0$ when $u\in\mathcal{F}(Q_T)$ satisfies \eqref{e5.3} with $g=0$.  We introduce the solution $v$ of the adjoint system as in \eqref{e5.7}.  Then $v\in\mathcal{Y}_2(Q_T)$ and
\begin{equation*}
\begin{split}
&\int_{Q_{T}}|u|^2dxdt =\int_{Q_{T}} u\enspace (-\displaystyle\frac{\partial v}{\partial t}-\Delta v+\nabla q)\enspace dxdt\\
&=(\text{ using \eqref{e5.13} with }\gamma_1(u)=0)\\
&=\int_{Q_{T}}u\enspace \nabla q\enspace dxdt =\int^T_0\int_\Omega u\enspace\nabla q\enspace dxdt,
\end{split}
\end{equation*}
and this last expression vanishes as in \eqref{e3.2a}.

The theorem is proved.

%
%

\section*{Acknowledgments.}  This work was partially
supported by the National Science Foundation under the grant
NSF-DMS1206438 and NSF-DMS1510249, and by the Research Fund of Indiana University.

\par


\begin{thebibliography}{33}

\bibitem[BDT13]{BDT13}
C. Bardos, F. Di Plinio, and R. Temam, The Euler equations in planar nonsmooth convex domains,
\textit{J. Math. Anal. and Applications}, \textbf{407}, 2013, 69--89.

\bibitem[BD73]{BD73}
J.D. Bozeman and C. Dalton, Numerical Study of viscous flow in a cavity,
\textit{Journal of Computational Physics}, Vol. \textbf{12}, 1973, 348--363.
\smallskip

\bibitem[Bot12]{Bot12}
Olivier Botella, Numerical solution of Navier-Stokes singular problem by a Chebyshev projection method, Thesis, Universit\'e de Nice, France, 2012.
\smallskip

\bibitem[Bur66]{Bur66}
O.R. Burgaff, Analytical and Numerical studies of the structures of steady separated flows, \textit{JFM}, Vol. \textbf{24}, 1966, 113--151.
\smallskip

\bibitem[Cat61]{Cat61}
Lamberto Cattabriga,
Su un problema al contorno relativo al sistema di equazioni di Stokes,
\textit{Rendiconti del Seminario Matematico della Universit\`a di Padova. The Mathematical Journal of the University of Padova},
\textbf{31}, 1961, 308--340.
\smallskip

\bibitem[CQT11a]{CQT11a} Q. Chen, Z. Qin and R. Temam, Numerical resolution near $t=0$ of nonlinear evolution equations in the presence of corner
singularities in space dimension 1, \textit{Communications in Computational Physics, CiCP}, \textbf{9}, No. 3, 2011, 568--586.
\smallskip

\bibitem[CQT11b]{CQT11b}
Q. Chen, Z. Qin and R. Temam, Treatment of incompatible initial and boundary data for parabolic equations in higher dimension,
\textit{Mathematics of Computation}, \textbf{80}, 276, 2011, 2071--2096.
\smallskip

\bibitem[Cuv86]{Cuv86}
C. Cuvelier, A. Segal and A.A. van Steenhoven, \textit{Finite element methods and Navier-Stokes equations}, Mathematics and its Applications, \textbf{22},
D. Reidel Publishing Co., Dordrecht, 1986, xvi+483.
\smallskip

\bibitem[Dav06]{Dav06}
T. Davis, \textit{Direct methods for sparse linear systems}, Fundamental of Algorithms, \textbf{2}, Society for Industrial and Applied Mathematics (SIAM), Philadelphia, PA, 2006, xii+217.
\smallskip

\bibitem[DL54]{DL54}
J. Deny and J.L. Lions,
Les espaces du type de Beppo Levi,
\textit{Universit\'e de Grenoble. Annales de l'Institut Fourier},
\textbf{5}, 1953--54, 305--370.
\smallskip

\bibitem[DT15]{DT15}
F. Di Plinio and R. Temam, Grisvard's shift theorem near $L^\infty$ and Yudovitch theory in polygonal domains, \textit{SIAM J. of Nonlinear Analysis},
\textbf{47}, No. 1, 2015, 159--178. DOI: 10.1137/130942632

\bibitem[Dun85]{Dun85}
D.A. Dunavant, High degree efficient symmetrical Gaussian quadrature rules for the triangle, \textit{International Journal for Numerical Methods in Engineering}, \textbf{21}, 6, 1985, 1129--1148.
\smallskip

\bibitem[ECG05]{ECG05}
E. Erturk, T.C. Corke, and C. Gokcol, Numerical Solutions of 2-D Steady Incompressible Driven Cavity Flow at High Reynolds Numbers,
\textit{International Journal for Numerical Methods in Fluids}, Vol. \textbf{48}, 2005, 747--774.
\smallskip

\bibitem[ESW05]{ESW05}
H. Elman, D. Silvester and A. Wathen, \textit{Finite elements and fast iterative solvers, with application in incompressible fluid dynamics},
Oxford University Press, Oxford, 2005.
\smallskip

\bibitem[FPT71a]{FPT71a}
M. Fortin, R. Peyret, and R. Temam, Calcul des \'ecoulements d'un fluide visqueux incompressible,
\textit{Proc. Second Intern. Conf. on Numerical Methods in Fluids Dynamics},
Lecture Notes in Physics, vol. \textbf{8}, Springer-Verlag, 1971.
\smallskip

\bibitem[FPT71b]{FPT71b}
M. Fortin, R. Peyret and R. Temam,
R\'esolution num\'erique des \'equations de Navier-Stokes pour un fluide incompressible,
\textit{J. M\'ecanique}, \textbf{10}, 1971, 357--390.
\smallskip

\bibitem[FQT12]{FQT12} N. Flyer, Z. Qin and R. Temam, A penalty method for numerically handling dispersive equations with incompatible initial
and boundary data, \textit{Numerical Methods for Partial Differential Equations}, \textbf{28}, 6, 2012, 1996-2009.  DOI:  10.1002/num.21693
\smallskip

\bibitem[Gar11]{Gar11}
Salvador Garcia, Aperiodic, chaotic lid-driven square cavity flows,
\textit{Journal Mathematics and Computers in Simulation}, \textbf{81} 9, 2011, 1741--1769.
Elsevier Science Publishers B. V. Amsterdam, The Netherlands, The Netherlands
table of contents doi>10.1016/j.matcom.2011.01.011
\smallskip

\bibitem[GGS82]{GGS82}
U. Ghia, K.N. Ghia, and C.T. Shin,
High-Resolutions for incompressible flow using the Navier-Stokes equations and a multigrid method,
\textit{Journal of Computational Physics}, Vol. \textbf{48}, 1982, 387--411.

\bibitem[Ghi84]{Ghi84}
J.M. Ghidaglia, R\'egularit\'e des solutions de certains probl\`emes aux limites lin\'eaires li\'es aux \'equations d'Euler,
\textit{Comm. in Partial Differntial Equations}, \textbf{9}, (13), 1984, 1265--1298.
\smallskip

\bibitem[Gri92]{Gri92}
P. Grisvard,
Singularities in boundary value problem,
Recherches en Math\'ematiques Appliqu\'ees,
Masson, Paris, Springer-Verlag, Berlin,
\textbf{22}, 1992, xiv+199
\smallskip

\bibitem[Gri85]{Gri85}
P. Grisvard,
\textit{Elliptic problems in nonsmooth domains},
Monographs and Studies in Mathematics,
Pitman (Advanced Publishing Program), Boston, MA,
\textbf{24}, 1985, xiv+410.
\smallskip

\bibitem[GMPQ02]{GMPQ02}
J.-L. Guermond, C. Migeon, G. Pineau, and L. Quartapelle,
Start-up flows in a three-dimensional rectangular driven cavity of aspect ratio
1:1:2 at Re$=1000$, \textit{Journal of Fluid Mechanics},
\textbf{450}, 2002, 169--199.
\smallskip

\bibitem[GMN81]{GMN81}
Murli M. Gupta, Ram P. Manohar, and Ben Noble,
Nature of viscous flows near sharp corners,
\textit{Computers and Fluids}, \textbf{9}, 4, 1981, 379--388.
\smallskip

\bibitem[HS15]{HS15}
Makram Hamouda, and Abir Sboui,
Boundary layers generated by singularities in the source function,
\textit{Asymptotic Analysis},  \textbf{93} (2015), no. 4, 281-310.​
\smallskip

\bibitem[Her81]{Her81}
B. H\'eron,
Quelques propri\'et\'es des applications de trace dans des espaces de champs de vecteurs \`a divergence nulle,
\textit{Communications in Partial Differential Equations}, \textbf{6}, 12, 1981, 1301--1334. DOI: 10.1080/03605308108820212
\smallskip

\bibitem[Kel83]{Kel83}
M.A. Kelmanson,
Modified integral equation solution of viscous flows near sharp corners, \textit{Computers \& Fluids}, \textbf{11}, 4, 1983, 307--324.
\smallskip

\bibitem[LM72]{LM72}
J.L. Lions, and E. Magenes,
\textit{Non-homogeneous boundary value problems and applications, Vol. I},
Springer-Verlag, New York-Heidelberg, 1972, xvi+357.
\smallskip

\bibitem[Mof63]{Mof63}
H.K. Moffatt, Viscous and resistive eddies near a sharp corner,
\textit{Journal of Fluid Mechanics}, \textbf{18}, 1, 1963, 1--18.
\smallskip

\bibitem[QT11b]{QT11b}
Q. Chen, Z. Qin and R. Temam, Treatment of incompatible initial and boundary data for parabolic equations in higher dimension,
\textit{Mathematics of Computation}, \textbf{80}, 276, 2011, 2071-2096.
\smallskip

\bibitem[QT12]{QT12} Z. Qin and R. Temam, Penalty method for the KdV equations, \textit{Applicable Analysis}, \textbf{91}, No. 1-2, 2012, 193-211,
DOI:  10.1080/00036811.2011.579564.
\smallskip

\bibitem[Rha84]{Rha84}
G. de Rham, \textit{Differentiable manifolds}, Grundlehren der Mathematischen Wissenschaften [Fundamental Principles of Mathematical Sciences], \textbf{266}, Springer-Verlag, Berlin, 1984, x+167
\smallskip

\bibitem[SLB89]{SLB89} W.W. Schultz, N.Y. Lee, N.Y. and J.P. Boyd,
Chebyshev pseudospectral method of viscous flows with corner singularities, \textit{Journal of Scientific Computing},
\textbf{4}, 1, 1989, 1--23.
\smallskip

\bibitem[SSB91]{SSB91}
M.R. Schumak, and W.W. Schultz,
Spectral method solution of the Stokes equations on nonstaggered grids, \textit{Journal of Computational Physics}, \textbf{94}, 1991, 30--58.
\smallskip

\bibitem[Ser83]{Ser83}
Denis Serre,
Equations de Navier-Stokes stationnaires avec donn\'ees peu r\'eguli\'eres,
\textit{Annali della Scuola Normale Superiore di Pisa. Classe di Scienze. Serie IV},
\textbf{10}, 4, 1983, 543--559.
\smallskip

\bibitem[She89]{She89}
Jie Shen,
\textit{Numerical simulation of the regularized driven cavity flows at high Reynolds numbers},
Spectral and high order methods for partial differential equations (Como, 1989),
273--280, North-Holland, Amsterdam, 1990
\smallskip

\bibitem[Sol64a]{Sol64a}
V.A. Solonnikov, Estimates for solutions of a non-stationary linearized system of Navier-Stokes equations, \textit{Akademiya Nauk SSSR. Trudy Matematicheskogo Instituta imeni V.A. Stekova}, \textbf{70}, 1964, 213--317.
\smallskip

\bibitem[Sol64b]{Sol64b}
V.A. Solonnikov, On the differential properties of the solution of the first boundary-value problem for a non-stationary system of Navier-Stokes equations, \textit{Akademiya Nauk SSSR. Trudy Matematicheskogo Instituta imeni V.A. Stekova}, \textbf{73}, 1964, 222--291.

\smallskip

\bibitem[Tem01]{Tem01}
Roger Temam,
\textbf{Navier-Stokes Equations Theory and Numerical Analysis},
North-Holland Pub. Company, in English, 1977, 500 pages.
Second revised edition, 1979. Third revised edition, 1984.
\par\noindent Russian Translation, Mir, Moscow, 1981.
Reedition in the AMS-Chelsea Series, AMS,Providence, 2001.
\smallskip

\bibitem[Wan93]{Wan93}
Xiaoming Wang, A remark on the characterization of the gradient of a distribution, \textit{Applicable Analysis}, \textbf{51}, 1993, 35--40.
\smallskip

\bibitem[Zha15]{Zha15}
Le Zhang, PhD in preparation, Indiana University, 2015.
\smallskip


\end{thebibliography}
\end{document}